\documentclass[leqno,12pt]{article}
\usepackage{amsthm,amsmath,amssymb}

\theoremstyle{plain}
\newtheorem{thm}{Theorem}
\newtheorem{cor}{Corollary}

\theoremstyle{definition}
\newtheorem{rem}{Remark}

\theoremstyle{remark}

\def\Cset{\mathbb{C}}

\title{A NEW SEMILOCAL CONVERGENCE THEOREM FOR THE WEIERSTRASS METHOD FROM DATA AT ONE POINT\thanks{This paper is published in: C. R. Acad. Bulg. Sci 59 (2006), No 2, 131--136.}}
\author{Petko~D.~Proinov}
\date{}

\begin{document}
\maketitle

\begin{abstract}
In this paper we present a new semilocal convergence theorem from data at one point for the Weierstrass iterative method for the simultaneous computation of polynomial zeros. The main result generalizes and improves all previous ones in this area.


\end{abstract}

\section{Introduction}

Let $f$  be a monic polynomial of degree $n \ge 2$ with simple complex zeros.
We consider the roots of $f$ as a point in ${\Cset}^n$.
Namely, a point $\xi$ in ${\Cset}^n$  with distinct coordinates is said to be a
root-vector of $f$ if each of its coordinates is a zero of $f$.
Starting from an initial point $z^0$ in ${\Cset}^n$ with distinct coordinates
we build in ${\Cset}^n$ the \textsc{Weierstrass} iterative sequence \cite{Wei03}
\begin{equation}  \label{Wei}
z^{k + 1}  = z^k - W(z^k ), \qquad k = 0,1,2, \ldots,
\end{equation}
where the operator $W$ in ${\Cset}^n$ is defined by
$W(z) = (W_1(z), \ldots, W_n(z))$ with
\[
W_i(z) = \frac{f(z_i)}{\prod \limits_{j \ne i}{(z_i - z_j)}}
\qquad (i = 1, 2, \cdots, n).
\]

It is well-known that under some initial conditions the Weierstrass sequence
\eqref{Wei} is well-defined and tends to a root-vector of $f$.
Iteration formula \eqref{Wei} defines the famous Weierstrass method
(known also as the Durand-Dochev-Kerner-Pre\v si\'c method) for finding all the zeros of $f$ simultaneously.

In 1962, \textsc{Dochev} \cite{Doc62,ID63} proved the first local convergence theorem for the Weierstrass method. Since 1980 a number of authors
\cite{Pre80,Zhe82,Zhe87,ZW93,WZ95,PCT95,Pet96,PH96,PHI98,Bat98,Han00,PH01}
have obtained semilocal convergence theorems for the Weierstrass method from data at one point (point estimation).
In this note we present a new semilocal convergence theorem for the Weierstrass method which improves and generalizes all these results.

The main result (Theorem~\ref{thm:Wei}) of this note will be proved elsewhere.

\section{Statement of the main result}

Throughout the paper the norm $\|.\|_p$ in ${\Cset}^n$
is defined as usual, i.e.
$\|z\|_p = \left( \sum\nolimits_{i = 1}^n {|z_i|^p} \right)^{1/p}$. For a given point $z$ in ${\Cset}^n$ with distinct coordinates we use the notations
\[
\frac{W(z)}{d(z)} =
\left( {\frac{W_1(z)}{d_1(z)}, \ldots, \frac{W_n(z)}{d_n(z)}} \right)
\quad\text{and}\quad E(z) = \left\|  \frac{W(z)}{d(z)} \right\|_p,
\]
where $d(z) = (d_1(z), \ldots, d_n(z))$ and
\[
d_i(z) = \min_{j \ne i} |z_i  - z_j | \qquad (i = 1,2, \ldots, n).
\]

\begin{thm} \label{thm:Wei} 
Let $f$ be a monic polynomial of degree $n \ge 2$ with simple zeros.
Let $1 \le p \le \infty$ and
$1/p + 1/q = 1$.
Define the real function
\begin{equation} \label{eq:phiW}
\phi(x) = \frac{(n - 1)^{1/q} x}{(1 - x)(1 - 2^{1/q} x)}
\left (1 + \frac{x}{(n - 1)^{1/p} (1 - 2^{1/q} x)} \right )^{n - 1} .
\end{equation}
Suppose that $z^0$ is an initial point in ${\Cset}^n$ with distinct coordinates  satisfying
\begin{equation} \label{iniW}
E(z^0) < 1/2^{1/q} \quad\text{and}\quad \phi(E(z^0)) \le 1 .
\end{equation}
Then the following statements hold true.
\begin{enumerate}
\item The Weierstrass iterative sequence \eqref{Wei} is well-defined and convergent to a root-vector $\xi $ of $f$. Moreover, the convergence is quadratic if
$\phi(E(z^0 )) < 1$.
\item For each $k \ge 1$ we have the following a priori error estimate
\begin{equation}  \label{priW}
\left \| z^k - \xi \right \|_p  \le \frac{{\theta} ^k {\lambda }^{2^k - 1}}
{1 - \theta {\lambda} ^{2^k } } \left \| z^1  - z^0 \right \|_p ,
\end{equation}
where $\lambda = \phi(E(z^0))$ and $\theta  = 1 - 2^{1/q} E(z^0 )$.
\item For all $k \ge 0$ we have the following a posteriori error estimate
\begin{equation}  \label{posW}
\left \| z^{k + 1} - \xi \right \|_p  \le \frac{{\theta}_k {\lambda}_k}
{1 - {\theta}_k {\lambda}_k^2} \left \| z^{k + 1} - z^k \right \|_p ,
\end{equation}
where $\lambda_k  = \phi(E(z^k))$ and $\theta_k = 1 - 2^{1/q} E(z^k)$.
\end{enumerate}
\end{thm}

\begin{rem}
Let $R(n,p)$ denote the unique solution of the equation $\phi(x) = 1$ in the interval $(0, 1/2^{1/q})$, where $\phi$ is defined by \eqref{eq:phiW}. Then the assumption \eqref{iniW} of Theorem~1 can also be written in the form
$E(z^0) \le R(n,p)$.
\end{rem}

\section{Comparison with the previous results}

In this section we compare Theorem~\ref{thm:Wei} with all previous results of the same type.
Note that for $z \in \Cset^n$ with distinct coordinates we have the obvious inequality
\begin{equation}  \label{dd}
\left \| \frac{W(z)}{d(z)} \right \|_p \le \frac{\| {W(z)} \|_p}{\delta (z)}
\quad\text{where}\quad \delta(z) = \min \{ d_1(z), \ldots, d_n(z) \}.
\end{equation}

\begin{cor} \label{cor:WeiALp}  
Let $f$ be a monic polynomial of degree $n \ge 2$ with simple zeros.
Let $1 \le p \le \infty$ and $1/p + 1/q = 1$.
Suppose $z^0$ is an initial point in
$\Cset^n $ satisfying
\begin{equation}  \label{ini1WLp}
\left \| \frac{W(z^0)}{d(z^0)} \right \|_p \le
\frac{1}{A(n - 1)^{1/q} + 2^{1/q} + 1} \, ,
\end{equation}
where $A = 1.763222 \ldots$  is the unique solution of the equation
$\exp{(1/x)} = x$.
Then the Weierstrass sequence \eqref{Wei} is quadratically convergent to a
root-vector $\xi$ of $f$. Moreover, we have the error estimates \eqref{priW} and \eqref{posW}.
\end{cor}

\begin{proof}
The statement of Corollary~\ref{cor:WeiLp} holds true even if we replace the
right-hand side of \eqref{ini1WLp} by
\(
R = 2 / ( m + \sqrt{m^2 - 4b} \, ),
\)
where $m = A a + b + 1$, $a = (n - 1)^{1/q}$ and $b = 2^{1/q}$. Define
\[
g(x) = \frac{ax}{(1 - x)(1 - bx)}.
\]
It is easy to show that  $g(R) = 1/A$ and $\phi (x) < g(x) \exp{g(x)}$ for
$0 < x < 1/b$.
Then by the definition of $A$ we get $\phi (R) < 1$ which according to
Theorem~\ref{thm:Wei} completes the proof.
\end{proof}

\textsc{Batra} \cite{Bat98} has proved that the Weierstrass method is convergent under the condition
\(
\| W(z^0) \|_{\infty} < \delta(z^0) / (2n).
\)
The following corollary improves Batra's result as well as previous
results \cite {Pre80,PCT95,Pet96,PH96,PHI98}.

\begin{cor} \label{cor:WeiLp}  
Let $f$ be a monic polynomial of degree $n \ge 2$ with simple zeros.
Let $1 \le p \le \infty$ and $1/p + 1/q = 1$.
Suppose $z^0$ is an initial point in
$\Cset^n $ satisfying
\begin{equation}  \label{ini2WLp}
\left \| \frac{W(z^0)}{d(z^0)} \right \|_p \le
\frac{1}{2(n - 1)^{1/q} + 2} .
\end{equation}
Then the Weierstrass sequence \eqref{Wei} is convergent to a
root-vector $\xi$ of $f$. Moreover, we have the error estimates \eqref{priW} and \eqref{posW}.
\end{cor}

\begin{proof}
Denote by $R$  the right-hand side of \eqref{ini2WLp}, i.e. $R = 1/(2a + 2)$, where again $a = (n - 1)^{1/q}$.
By Theorem~\ref{thm:Wei} it suffices to prove that $\phi(R) \le 1$. If $n = 2$, then
$\phi(R) = (4/3)(5-b) / (4-b)^2 \le 1$,
where $b=2^{1/q}$ . Further, suppose $n \ge 3$. It is easy to see that
\begin{equation}  \label{phiR}
\phi(R) \le \frac{2a(a + 1)}{(2a + 1)(2a + 2 - b)} \sqrt e .
\end{equation}
If  $a \ge 2$, then \eqref{phiR} implies
\[
\phi(R) \le \frac{a + 1}{2a + 1} \sqrt e  \le \frac{3}{5} \sqrt e  < 1.
\]
If  $a \le 2$, then \eqref{phiR} implies
\[
\phi(R) \le \frac{2a(a + 1)}{(2a + 1)(a + 2)} \sqrt e \le
\frac{3}{5} \sqrt e < 1
\]
which completes the proof.
\end{proof}

\begin{cor} \label{cor:WeiL1}  
Let $f$ be a monic polynomial of degree $n \ge 2$ with simple zeros.
Suppose that $z^0$ is an initial point in $\Cset^n $ satisfying
\[
\left \| \frac{W(z^0)}{d(z^0)} \right \|_{1} \le R = 0.307541 \ldots ,
\]
where $R$  is the unique solution of the equation
\begin{equation}  \label{iniWL1}
{\left( \frac{x}{1-x} \right)}^2 \exp{\frac{x}{1-x}} = 1.
\end{equation}
in the interval $(0,1)$.
Then the conclusion of Corollary~\ref{cor:WeiALp} holds for $p=1$.
\end{cor}

\begin{proof}
It is easy to show that $\phi(x) < g(x)$ for $0 < x < 1$, where $g(x)$ denotes
the left-hand side of equation \eqref{iniWL1}.
Therefore, $\phi(R)< 1$ which according to Theorem~\ref{thm:Wei} completes the proof.
\end{proof}

\begin{cor}[Han \cite{Han00}] \label{cor:WeiHan}  
Let $f$ be a monic polynomial of degree $n \ge 2$ with simple zeros.
Assume $1 \le p \le \infty$ and $1/p + 1/q = 1$.
Suppose that
\begin{equation}  \label{iniHan}
\left \| \frac{W(z^0)}{d(z^0)} \right \|_p \le \frac{n(2^{1/n} - 1)}{(n - 1)^{1/q}  + 2^{1/q}}
\left (1 - \frac{(n - 1)^{1/q} \, n(2^{1/n} - 1)}{(n - 1)^{1/q} + 2^{1/q}} \right).
\end{equation}
Then the Weierstrass iterative sequence \eqref{Wei} converges to a root-vector of $f$.
\end{cor}

\begin{proof}
Let $R$ denote the right-hand side of \eqref{iniHan}. According to
Corollary~\ref{cor:WeiLp} it suffices to prove the inequality
$R \le 1/(2a + 2)$,
where $a$ and $b$ are defined as in the proof of Corollary~\ref{cor:WeiALp}.
Write $R$ in the form $R = g(\tau)$, where $g(x) = x(1 - a x)$ and
$\tau  = n(2^{1/n} - 1) / (a + b)$.
Since $n(2^{1/n} - 1) < 1$ we get $R \le g(1/(2a)) = 1/(4a) \le 1/(2a+2)$ and the proof of the corollary is complete.
\end{proof}

Setting $p = \infty$ in Theorem~\ref{thm:Wei} and taking into account \eqref{dd} we obtain the following corollary. The first part of it is due to
\textsc{Zheng} \cite{Zhe82}
and \textsc{Petkovi\'c} and \textsc{Herceg} \cite{PH01}. The second part of the corollary is due to
\textsc{Zheng} \cite{Zhe87}.

\begin{cor}[Zheng \cite{Zhe82,Zhe87} and Petkovi\'c and Herceg \cite{PH01}] \label{cor:WeiZPH}  
Let $f$ be a monic polynomial of degree $n \ge 2$ with simple zeros.
Let $0<C< \frac{1}{2}$ and
\[
\lambda : = \frac{(n - 1)C}{(1 - C)(1 - 2C)}
\left(1 + \frac{C}{1 - 2C} \right)^{n - 1} \le 1
\quad\text{where}\quad C =  \frac{\| W(z^0) \|_{\infty}}{\delta(z^0)}.
\]
Then the Weierstrass method \eqref{Wei} is convergent to a root-vector $\xi$ of
$f$. Moreover, the error estimate \eqref{priW} holds with $\theta  = 1 - 2C$ and $p=\infty$.
\end{cor}

Petkovi\'c and Herceg \cite{PH01} have proved that the Weierstrass method is convergent under the condition
\(
\| W(z^0) \|_{\infty} < \delta(z^0) / (an+b),
\)
where $a=1.76325$ and $b=0.8689425$.
Note that for sufficiently large $n$ ($n \ge 13624$) this result is an immediate consequence of Corollary~\ref{cor:WeiALp}.
The following corollary improves Petkovi\'c and Herceg's result for all
$n \ge 2$.

\begin{cor} \label{cor:WeiInf}  
Let $f$ be a monic polynomial of degree $n \ge 2$ with simple zeros. Under the initial condition
\[
\left \| \frac{W(z^0)}{d(z^0)} \right \|_{\infty} \le
\frac{1}{1.76325 n + 0.6869}
\]
the Weierstrass method is convergent to a root-vector $\xi$ of $f$ with the second order of convergence. Moreover, the estimates \eqref{priW} and
\eqref{posW} hold for $p=\infty$.
\end{cor}

\begin{proof}
The sequence $\phi_n = \phi(1/(1.76325 n + 0.6869))$ is increasing for $2 \le n \le 132$ and decreasing for $n \ge 132$. Hence $\phi_n \le \phi_{132} < 1$ for
$n \ge 2$. Now the conclusion follows from Theorem~\ref{thm:Wei}.
\end{proof}

\begin{cor}[Wang and Zhao \cite{WZ95}] \label{cor:WeiInfWZ} 
Let $f$ be a monic polynomial of degree $n \ge 2$ with simple zeros. Let us assume that
\begin{equation}  \label{Winfd}
\| W(z^0) \|_{\infty} \le C(n) \delta(z^0)
\quad\text{where}\quad
C(n) = - \mathop{\min} \limits_{x > 0} (x(1 + x)^{n - 1} - 2x).
\end{equation}
Then the Weierstrass iterative sequence \eqref{Wei} converges to a root-vector of $f$.
\end{cor}

\begin{proof}
By Corollary~\ref{cor:WeiLp} and \eqref{dd} it suffices to prove that $C(n) \le 1/(2n)$. This is obvious for $n = 2$ and $n = 3$ since $C(2) = 0.25$ and
$C(3) = 0.112 \ldots$ Define $g(x) = 2x - x(1 + x)^{n - 1}$. It is easy to show that there exists a unique real number $t$ such that $0 < t < 2^{1/(n - 1)} - 1$ and $g'(t) =0$.
One can prove that $t < 1/(2n - 1)$ for $n \ge 4$. Then by the definitions of
$C(n)$ and $t$ we get
\[
C(n) = \mathop {\max }\limits_{x \in \, [ \, 0 \, , \, 2^{1/(n - 1)} - 1]}
g(x) = g(t)
= 2(n - 1) t^2 / (1 + n t)
\]
which implies  $C(n) < 1 / (3n)$ for $n \ge 4$.
This completes the proof of Corollary~\ref{cor:WeiInfWZ}.
\end{proof}

\begin{cor}[Wang and Zhao \cite{WZ95}] \label{cor:WeiL1WZ} 
Let $f$ be a monic polynomial of degree $n \ge 4$ with simple zeros. Let us suppose that
\[
\| W(z^0) \|_{1} \le C(n) \; \delta(z^0)
\quad\text{where}\quad C(n) =  - \mathop{\min} \limits_{x > 0}
\left( \sum\nolimits_{j = 1}^{n - 1}{\frac{n - j}{j! n} x^{j + 1}} - x \right).
\]
Then the Weierstrass iterative sequence \eqref{Wei} converges to a root-vector of $f$.
\end{cor}

\begin{proof}
According to Corollary~\ref{cor:WeiL1} and \eqref{dd} it suffices to prove that
$C(n) \le 0.3$. For a given $n \ge 2$ define the function
\(
f_n(x) = \sum\nolimits_{j = 1}^{n - 1} {\frac{n - j}{j! n} x^{j + 1}} - x.
\)
It is easy to verify that for $x > 0$ we have
\[
f_n(x) \le \sum\nolimits_{j = 1}^n {\frac{x^{j + 1}}{j!}} - x
= (n + 1)f_{n + 1}(x) - nf_n(x).
\]
This implies $f_n(x) \le f_{n + 1}(x)$. Now by the definition of $C(n)$ we get
$C(n + 1) \le C(n)$ for $n \ge 2$. Hence $C(n) \le C(4) = 0.279 \ldots$ for $n \ge 4$ which completes the proof.
\end{proof}

\begin{rem}
\textsc{Wang} and \textsc{Zhao} \cite{ZW93} have also proved that the Weierstrass method
\eqref{Wei}
is convergent under the condition
\[
\left\| W(z^0) \right\|_{1} \le \frac{(3 - 2 \sqrt 2) n}{n-1} \; \delta(z^0 ).
\]
Note that for $n \ge 3$ this result is an immediate consequence of
Corollary~\ref{cor:WeiALp} and
\eqref{dd}.
\end{rem}

\section{Remark on the SOR Weierstrass method}
In this section we are concerned with the successive overrelaxation (SOR) Weierstrass method
\begin{equation}  \label{SORW}
z^{k + 1}  = z^k - h_k \, W(z^k ), \qquad k = 0,1,2, \ldots,
\end{equation}
where $h_k = h_k(f) \in (0,1]$ is an acceleration parameter.
\textsc{Petkovi\'c} and \textsc{Kjurkchiev} \cite{PK97} have noticed that usually the SOR method
\eqref{SORW} is faster if $h_k$ is closer to 1.
In 1995, \textsc{Wang} and \textsc{Zhao} \cite{WZ95} considered the SOR method \eqref{SORW} with $h_k$ defined by
\begin{equation}  \label{hkWZ}
 h_k  = \min \left\{ {1,\frac{{0.204378 \:
\delta(z^k)}}{{\sum\limits_{i = 1}^n {|W_i (z^k )|} }}} \right\},
\end{equation}
According to Corollary~\ref{cor:WeiL1} we can consider the SOR Weierstrass method \eqref{SORW} with $h_k$ defined by
\begin{equation}  \label{hkP}
h_k  = \min \left\{ 1,\frac{0.307541}{\sum\limits_{i = 1}^n
{\left| \frac{W_i(z^k )}{d_i(z^k )} \right|}}  \right\}.
\end{equation}
Note that the new $h_k$ is closer to 1. Moreover, if our acceleration parameter $h_k$ is less than 1, then it is greater than Wang-Zheng's parameter $h_k$ by more than 50\%.

\noindent
Faculty of Mathematics and Informatics\\
University of Plovdiv\\
Plovdiv 4000, Bulgaria\\
E-mail: proinov@pu.acad.bg


\begin{thebibliography}{99}

\bibitem{Wei03}
\textsc{K.~Weierstrass},
Neuer Beweis des Satzes, dass jede ganze rationale Funktion einer Veranderlichen dargestellt werden kann als ein Produkt aus linnearen Functionen derselben Veranderlichen,
Ges. Werke \textbf{3} (1903), 251--269.

\bibitem{Doc62}
\textsc{K.~Dochev},
Modified Newton method for simultaneous approximation of all roots of a given algebraic equation,
Phys. Math. J. Bulg. Acad. Sci. \textbf{5} (1962), 136--139 (in Bulgarian).

\bibitem{ID63}
\textsc{L.~Ilieff} and \textsc{K.~Dochev},
{\"U}ber Newtonsche Iterationen.
Wiss. Z. Tech. Univ. Drezden \textbf{12} (1963), 117--118.

\bibitem{Pre80}
\textsc{M.~D.~Pre\v si\'c},
A convergence theorem for a method for simultaneous determination of all zeros of a polynomial,
Publ. Inst. Math. (N.S.) \textbf{28} (1980), 159--165.

\bibitem{Zhe82}
\textsc{S.~M.~Zheng},
On convergence of the Durand-Kerner's method for finding all roots of a polynomial simultaneously,
Kexue Tongbao \textbf{27} (1982), 1262--1265.

\bibitem{Zhe87}
\textsc{S.~M.~Zheng},
On convergence of a parallel algorithm for finding the roots of a polynomial,
J. Math. Res. Exp. \textbf{7} (1987), 657--660 (in Chinese).

\bibitem{ZW93}
\textsc{F.~G.~Zhao} and \textsc{D.~R.~Wang},
The theory of Smale's point estimation and the convergence of Durand-Kerner
program,
Math. Numer. Sinica \textbf{15} (1993), 196--206 (in Chinese).

\bibitem{WZ95}
\textsc{D.~R.~Wang} and {F.~G.~Zhao},
The theory of Smale's point estimation and its applications,
J. Comput. Appl. Math. \textbf{60} (1995), 253--269.

\bibitem{PCT95}
\textsc{M.~Petkovi\'c}, \textsc{C.~Carstensen} and \textsc{M.~Trajkovi\'c},
Weierstrass formula and zero-finding methods,
Numer. Math. \textbf{69} (1995), 353--372.

\bibitem{Pet96}
\textsc{M.~S. Petkovi\'c},
On initial conditions for the convergence of simultaneous root-finding methods,
Computing \textbf{57} (1996), 163--177.

\bibitem{PH96}
\textsc{M.~S.~Petkovi\'c, D.~Herceg},
Point estimation and safe convergence of root-finding simultaneous methods,
Sci. Rev. \textbf{21-22} (1996), 117--130.

\bibitem{PHI98}
\textsc{M.~S.~Petkovi\'c, D.~Herceg, S.~Ili\'c},
Safe convergence of simultaneous methods for polynomial zeros,
Numer. Algorithms \textbf{17} (1998), 313--331.

\bibitem{Bat98}
\textsc{P.~Batra},
Improvement of a convergence condition for Durand-Kerner iteration,
J. Comput. Appl. Math. \textbf{96} (1998), 117--125.

\bibitem{Han00}
\textsc{D.~F.~Han},
The convergence of the Durand-Kerner method for simultaneously
finding all zeros of a polynomial,
J. Comput. Math. \textbf{18} (2000), 567--570.

\bibitem{PH01}
\textsc{M.~S. Petkovi\'c, D.~Herceg},
Point estimation of simultaneous methods for solving polynomial equations,
J. Comput. Appl. Math. \textbf{136} (2001), 283--307.

\bibitem{PK97}
\textsc{M.~S.~Petkovi\'c, N.~Kjurkchiev},
A note on the convergence of the Weierstrass SOR method for polynomial roots,
J. Comput. Appl. Math. \textbf{80} (1997), 163--168.

\end{thebibliography}
\end{document}